\documentclass{article}

\usepackage{amsmath}
\usepackage{amssymb}
\usepackage{graphicx}
\usepackage{setspace}
\usepackage{citesort}
\usepackage[squaren]{SIunits}
\usepackage{url}

\newcommand{\E}{\mathrm{e}}
\newcommand{\D}{\mathrm{d}}
\newcommand{\un}{u^{(n)}}
\newcommand{\fn}{f^{(n)}}

\begin{document}

\title{Numerical solution of the modified Bessel equation}

\author{Michael Carley\thanks{Department of Mechanical
    Engineering, University of Bath, Bath BA2 7AY, United Kingdom
    ({\tt{m.j.carley@bath.ac.uk})}}}

\bibliographystyle{plain}

\maketitle

\begin{abstract}
  A Green's function based solver for the modified Bessel equation has
  been developed with the primary motivation of solving the Poisson
  equation in cylindrical geometries. The method is implemented using
  a Discrete Hankel Transform and a Green's function based on the
  modified Bessel functions of the first and second kind. The
  computation of these Bessel functions has been implemented to avoid
  scaling problems due to their exponential and singular behavior,
  allowing the method to be used for large order problems, as would
  arise in solving the Poisson equation with a dense azimuthal
  grid. The method has been tested on monotonically decaying and
  oscillatory inputs, checking for errors due to interpolation and/or
  aliasing. The error has been found to reach machine precision and to
  have computational time linearly proportional to the number of
  nodes.
\end{abstract}


  

\pagestyle{myheadings}

\thispagestyle{plain}

\markboth{M. CARLEY}{Numerical solution of the modified Bessel equation}

\section{Introduction}
\label{sec:intro}


This paper is motivated by the requirement for a Poisson solver for
cylindrical domains. Such a solver is a basic element in solving a
range of physical problems and accurate methods have been developed
for Cartesian~\cite{genovese-deutsch-goedecker07} and
cylindrical~\cite{chen-su-shizgal00} domains. An issue which arises in
solving the problem in a cylindrical coordinate system is the
singularity which arises at the axis due to the form of the
differential operator. A recent paper by Pataki and
Greengard~\cite{pataki-greengard11} introduces a Green's function
solver for the Poisson equation, which avoids problems with this
singularity, and automatically imposes a radiation boundary condition,
by using an integral formulation for the solution. 

The method which Pataki and Greengard~\cite{pataki-greengard11}
develop is based on Fourier transforms in the axial and azimuthal
coordinate, followed by the solution of a modified Bessel equation in
the radial coordinate. They test the accuracy of their method on an
axisymmetric problem with monotonic decay in radius, and show its
application to an asymmetric problem. A problem which arises in their
algorithm is that the integration technique used does not work well
for large azimuthal orders, i.e.\ meshes dense in angle, or for large
axial wavenumbers, i.e.\ meshes dense in the axial coordinate, due to
the poor scaling of the modified Bessel functions which appear in the
Green's function for the problem.

This motivated an attempt to find a robust method for solving the
modified Bessel equation, which would work for large wavenumbers and
azimuthal orders. The method developed, which is described in the rest
of this paper, is based on the Discrete Hankel Transform (DHT),
tabulated integrals, and recursions for the modified Bessel functions.

\section{Problem formulation}
\label{sec:formulation}

The Poisson equation in cylindrical coordinates is:
\begin{align}
  \label{equ:poisson}
  u_{rr}(r,\theta,z) + \frac{1}{r}u_{r}(r,\theta,z) + 
  \frac{1}{r}u_{\theta\theta}(r,\theta,z) + 
  u_{zz}(r,\theta,z) &= f(r,\theta,z),
\end{align}
where $(r,\theta,z)$ are cylindrical coordinates, $u$ is the solution
and $f$ is some forcing term. With the problem defined on nodes
regularly spaced in $\theta$ and $z$, this equation can be
solved~\cite{pataki-greengard11} by using the FFT to Fourier transform
$u$ and $f$ in $\theta$ and $z$ to yield a set of modified Bessel
equations:
\begin{align}
  \label{equ:problem}
  \un_{rr}(r,\kappa) + \frac{1}{r}\un_{r}(r,\kappa) - 
  \left(
    \frac{n^{2}}{r^{2}} + \kappa^{2}
  \right)
  \un(r,\kappa) &= \fn(r,\kappa),
\end{align}
where $n$ is the azimuthal order and $\kappa$ the axial wavenumber,
with $\un$ and $\fn$ the Fourier transformed solution and forcing term
respectively. After solving Equation~\ref{equ:problem} for each value
of $n$ and $\kappa$, $\un(r,\kappa)$ can be inverse Fourier
transformed to yield the solution $u(r,\theta,z)$.

Pataki and Greengard~\cite{pataki-greengard11} give a method for
solving this modified Bessel equation, subject to a radiation boundary
condition at some outer radius $r=R$. This is done using the Green's
function for the modified Bessel equation with the solution written:
\begin{align}
  \label{equ:solution}
  \un(r,\kappa) &= \int_{0}^{R}G_{n}(\kappa, r, s)f(s, \kappa)\,\D s.
\end{align}
The Green's function $G_{n}$ is:
\begin{align}
  \label{equ:greens:func}
  G_{n}(\kappa, r, s) &=
  \left\{
    \begin{array}{ll}
      -s I_{n}(\kappa r)K_{n}(\kappa s),\quad &r\leq s;\\
      -s K_{n}(\kappa r)I_{n}(\kappa s), &r\geq s,
    \end{array}
  \right.
\end{align}
where $I_{n}(x)$ and $K_{n}(x)$ are the modified Bessel functions of
the first and second kind respectively. The first of these, $I_{n}$,
grows exponentially, while $K_{n}$ decays exponentially, but has
$x^{-n}\log x$ behavior at the origin, leading to some numerical
problems caused by the scaling of the Green's function. In practice,
the solution is computed as:
\begin{align}
  \un(r,\kappa) &= 
  -\frac{K_{n}(\kappa r)}{K_{n}(\kappa R_{m})}
  \int_{0}^{r}I_{n}(\kappa s)K_{n}(\kappa R_{m})\fn(s)s\,\D s \nonumber\\
  &-I_{n}(\kappa r)K_{n}(\kappa R_{m})
  \int_{r}^{R}\frac{K_{n}(\kappa s)}{K_{n}(\kappa R_{m})}\fn(s)s\,\D s,
  \label{equ:solution:scaled}
\end{align}
where the radial domain is divided at radii $R_{m}$ which are used to
set reference values of the modified Bessel functions. The product
$I_{n}(\kappa s)K_{n}(\kappa R_{m})$, and the ratio $K_{n}(\kappa
s)/K_{n}(\kappa R_{m})$, are thus well-scaled avoiding problems in
computation, as long as $|\kappa(s-R_{m})|$ is not too
large. Pataki~\cite{pataki11} reports that the integration is
performed using a dyadic grid, and that by scaling the modified Bessel
functions as they are computed, the method works well for
$n\lesssim40$, corresponding to~80 points in the azimuthal mesh.

For many applications, it is desirable to use a denser mesh than this
and so a different approach was sought for the solution of
Equation~\ref{equ:problem}. The natural transform technique for
problems in polar coordinates is the Discrete Hankel Transform (DHT),
which expands a function as a series of ordinary Bessel functions
$J_{m}(x)$. For a function $f(r)$, $0\leq r\leq R$:
\begin{align}
  \label{equ:dht}
  f(r) &\approx \sum_{m=1}^{M}\hat{f}_{m} J_{m}(\alpha_{m} r)
\end{align}
where $J_{m}(\alpha_{m} R)\equiv 0$ and $\hat{f}_{m}$ denotes the $m$
coefficient of the expansion of $f$. If the DHT of $\fn(r,\kappa)$ is
available, the solution of Equation~\ref{equ:problem} can be
immediately written:
\begin{align}
  \label{equ:dht:solution}
  \un(r,\kappa) &= \sum_{m=1}^{M}\hat{f}_{m} 
  \int_{0}^{R}J_{m}(\alpha_{m} s)G_{n}(\kappa, r, s)\,\D s.
\end{align}
The integrals required in Equation~\ref{equ:dht:solution} are given in
standard tables~\cite{gradshteyn-ryzhik80}:
\begin{subequations}
  \label{equ:quad:ij}
  \begin{align}
    \int_{0}^{r} I_{n}(\kappa s) J_{n}(\alpha s)s\,\D s
    &=
    \left[
      \alpha J_{n+1}(\alpha r)I_{n}(\kappa r) +
      \kappa I_{n+1}(\kappa r)J_{n}(\alpha r) 
    \right] \frac{r}{\alpha^{2}+\kappa^{2}},\\
    \int_{0}^{r} K_{n}(\kappa s) J_{n}(\alpha s)s\,\D s
    &=
    \left[
      \left(
        \frac{\alpha}{\kappa}
      \right)^{n}
      +
      \alpha r J_{n+1}(\alpha r)K_{n}(\kappa r) - 
      \kappa r K_{n+1}(\kappa r)J_{n}(\alpha r) 
    \right]\frac{1}{\alpha^{2}+\kappa^{2}} 
  \end{align}
\end{subequations}
which gives a solution for the problem in terms of the DHT
coefficients $\widehat{\fn}_{m}$ and $\kappa$. As written, this
solution is correct, but not numerically useful, due to the poor
scaling of the modified Bessel functions, especially for large values
of $\kappa$ and/or $n$. It must be rewritten in order to avoid
numerical difficulties. 

\section{Numerical implementation}
\label{sec:implementation}

In order to avoid numerical difficulties caused by poor scaling of the
modified Bessel function, Equations~\ref{equ:quad:ij} are rewritten
and used to give the convolution of the Green's function with the
ordinary Bessel function as:
\begin{subequations}
  \label{equ:numerical:IJK}
  \begin{align}
    \int_{0}^{R}&G_{n}(\kappa, r, s)J_{n}(\alpha s)\,\D s =\nonumber\\
    &-R\frac{I_{n}(\kappa r)K_{n}(\kappa R)}{\alpha^{2}+\kappa^{2}}
    \left[
      \alpha J_{n+1}(\alpha R) - 
      \kappa J_{n}(\alpha R)\frac{K_{n+1}(\kappa R)}{K_{n}(\kappa R)}
    \right]
    \nonumber\\
    &-\kappa r
    \frac{I_{n}(\kappa r)K_{n}(\kappa r)J_{n}(\alpha r)}{\alpha^{2}+\kappa^{2}}
    \left[
      \frac{I_{n+1}(\kappa r)}{I_{n}(\kappa r)}
      +
      \frac{K_{n+1}(\kappa r)}{K_{n}(\kappa r)}
    \right],\quad r\neq 0,\\
    &= \frac{1}{\alpha^{2}+\kappa^{2}},\quad r =0,\, n=0,\\
    &= 0, \quad r =0,\, n\neq 0.
  \end{align}
\end{subequations}
Written in this form, the modified Bessel functions appear only as
ratios $I_{n+1}(x)/I_{n}(x)$ and $K_{n+1}(x)/K_{n}(x)$ or as the
products $I_{n}(x)K_{n}(x)$ and $I_{n}(\kappa r)K_{n}(\kappa R)$. The
ratios can be computed directly using standard functional relations,
while the products are calculated using the same ratios combined with
modified Bessel functions of order zero, which can be computed
accurately and stably. Implementation of the solution technique thus
requires two elements, a method for the calculation of ratios of
modified Bessel functions, and a method for computing the DHT. In
practice, the input will not be defined on the nodes of the DHT, so an
interpolation scheme will also be required. The method has been coded
making use of the GNU Scientific Library
(GSL)~\cite{galassi-davies-etal05}, which provides functions for the
computation of scaled versions of the modified Bessel functions,
directly returning $I_{n}(x)\exp[-x]$ and $K_{n}(x)\exp[x]$. The
algorithm has been designed to use these scaled functions, to avoid
problems of underflow and overflow.

\subsection{Ratios of modified Bessel functions}
\label{sec:ratios}

The ratios of modified Bessel functions can be computed using standard
functional relations~\cite[8.486]{gradshteyn-ryzhik80}:
\begin{subequations}
  \label{equ:relations:IK}
  \begin{align}
    I_{n-1}(x) &= \frac{2n}{x}I_{n}(x) + I_{n+1}(x),\\
    K_{n+1}(x) &= \frac{2n}{x}K_{n}(x) + K_{n-1}(x),
  \end{align}
\end{subequations}
using an approach similar to that of Amos~\cite{amos74}, who writes
the ratios of successive functions as:
\begin{subequations}
  \label{equ:ratios:IK}
  \begin{align}
    \frac{I_{n}(x)}{I_{n-1}(x)} &= 
    \frac{x}{2n + x I_{n+1}(x)/I_{n}(x)},\\
    \frac{K_{n+1}(x)}{K_{n}(x)} &= 
    \frac{2n}{x} + \frac{K_{n-1}(x)}{K_{n}(x)}.
  \end{align}
\end{subequations}
The recursion for $I_{n}(x)/I_{n-1}(x)$ is stable for descending $n$
while that for $K_{n+1}(x)/K_{n}(x)$ is stable for increasing $n$. The
recursion for $K_{n+1}(x)/K_{n}(x)$ is seeded with
$K_{1}(x)/K_{0}(x)$, computed using the scaled form of $K_{0}(x)$ and
$K_{1}(x)$. The recursion for $I_{n}(x)/I_{n-1}(x)$ is seeded using
Olver's asymptotic formula~\cite[10.41.10]{dlmf10} for modified Bessel
functions of large order, as recommended by Amos~\cite{amos74},
starting at order equal to the larger of $n+8$ and~32. The asymptotic
expansion is given by:
\begin{align}
  \label{equ:olver}
  I_{n}(x) &\sim \frac{\E^{n\eta}}{(2\pi n)^{1/2}(1+z^{2})^{1/4}}
  \sum_{q=0}^{\infty}\frac{u_{q}(t)}{n^{q}},\\
  z &= x/n,\, \eta = (1+z^{2})^{1/2} + \log
  \frac{z}{1+(1+z^{2})^{1/2}},\,
  t = 1/(1+z^{2})^{-1/2},\nonumber\\
  u_{0}(t) &= 1, \nonumber\\
  u_{1}(t) &= (3t - 5t^{3})/24, \nonumber\\
  u_{2}(t) &= (81t^{2} - 462t^{4} + 385t^{6})/1152,\nonumber\\
  u_{3}(t) &= (30375t^{3} - 369603t^{5} + 765765t^{7} - 425425t^{9})/414720,
  \nonumber
\end{align}
while for small arguments, $(x/2)^{2}<n+1$, the series expansion of
$I_{n}(x)$ is used~\cite[8.445]{gradshteyn-ryzhik80}. 

Given a sequence of ratios of modified Bessel functions, the products
in Equation~\ref{equ:numerical:IJK} can be computed as:
\begin{align}
  \label{equ:product:IK}
  I_{n}(\kappa r)K_{n}(\kappa r) &=
  \left[
    I_{0}(\kappa r)\E^{-\kappa r}
  \right]
  \left[
    K_{0}(\kappa r)\E^{\kappa r}
  \right]
  \prod_{i=0}^{n-1}
  \left[
    \frac{I_{i+1}(\kappa r)}{I_{i}(\kappa r)}
  \right]
  \left[
    \frac{K_{i+1}(\kappa r)}{K_{i}(\kappa r)}
  \right],
\end{align}
and
\begin{align}
  \label{equ:product:IKR}
  I_{n}(\kappa r)K_{n}(\kappa R) &=
  \left[
    I_{0}(\kappa r)\E^{-\kappa r}
  \right]
  \left[
    K_{0}(\kappa R)\E^{\kappa R}
  \right]
  \prod_{i=0}^{n-1}
  A
  \left[
    \frac{I_{i+1}(\kappa r)}{I_{i}(\kappa r)}
  \right]
  \left[
    \frac{K_{i+1}(\kappa R)}{K_{i}(\kappa R)}
  \right],\\
  A &= \E^{\kappa(r-R)/n},\nonumber
\end{align}
where terms of the form $I_{n}(x)\exp[-x]$ and $K_{n}(x)\exp[x]$ are
computed directly using the scaled form of the modified Bessel
functions. The ratios of successive modified Bessel functions are well
scaled and multiplying them in pairs as in the products of
Equations~\ref{equ:product:IK} and~\ref{equ:product:IKR} avoids
underflow and overflow problems. 

\subsection{Discrete Hankel Transform}
\label{sec:dht}

The coefficients of the DHT are computed using the method of
Lemoine~\cite{lemoine94}. This is essentially a quadrature rule based
on the zeros of the ordinary Bessel function of order $n$,
$J_{n}(x)$. The function to be transformed is specified at these zeros
$x_{m}$, $0\leq m < M$, $J_{n}(x_{m})=0$, and the DHT is given by a
matrix multiplication of the vector of input data $f(x_{m})$ with the
matrix entries given by:
\begin{align}
  \label{equ:hankel}
  B_{mj}^{(n)} &= \frac{2}{x_{M}} 
  \frac{J_{n}(x_{m}x_{j}/x_{M})}{|J_{n+1}(x_{m})J_{n+1}(x_{j})|}.
\end{align}
In the calculations presented here, the GSL
implementation~\cite{galassi-davies-etal05} of Lemoine's method was
used, but with a modification to compute the zeros of $J_{n}(x)$ using
the $O(M)$ algorithm of Glaser et. al~\cite{glaser-liu-rokhlin07}. 

In order to compute the DHT, the input must be specified at the zeros
of the Bessel function. If only one order $n$ is of interest, this
presents no difficulties, but if the solution is to be found for
multiple values of $n$, as in solving a Poisson equation, for example,
an interpolation scheme is required to transfer the input from the
problem mesh onto the DHT nodes, in particular because the zeros are
not the same for different orders of Bessel function. 

\subsection{Interpolation}
\label{sec:interpolation}

Given that an interpolation scheme will almost always be needed, the
approach used by Pataki and Greengard~\cite{pataki-greengard11} has
been adopted. The domain $0\leq r\leq R$ is divided into $N$ blocks,
$R_{n}\leq r\leq R_{n+1}$, $n=0,\ldots,N-1$. Each block is discretized
with $P$ points, given by the Chebyshev nodes of the second
kind:
\begin{align}
  \label{equ:chebyshev}
  r_{nP+p} &= 
  \frac{R_{n+1}+R_{n}}{2} +
  \frac{R_{n+1}-R_{n}}{2}\cos\frac{p\pi}{P},\quad 
  p = 0,1,\ldots,P.
\end{align}
Evaluation of $\fn(r,\kappa)$ within each block is performed using
barycentric Lagrangian interpolation~\cite{berrut-trefethen04}. Since
Equation~\ref{equ:numerical:IJK} can be computed directly at arbitrary
values of $r$, the solution is generated on the input nodes, with no
requirement for interpolation from the DHT nodes.

\subsection{Summary of algorithm}
\label{sec:summary}

Given the elements described above, the solution algorithm can be
summarized as follows. For a given order $n$, wavenumber $\kappa$ and
input $\fn(r,\kappa)$, $0\leq r\leq R$:
\begin{enumerate}
\item generate, if necessary, the DHT matrix and corresponding nodes
  $r_{m}$;
\item if necessary, interpolate $\fn(r,\kappa)$ onto the nodes
  $r_{m}$;
\item perform the DHT to yield $\widehat{\fn}_{m}$;
\item evaluate Equation~\ref{equ:dht:solution} at the input nodes
  using Equations~\ref{equ:numerical:IJK}.
\end{enumerate}

\section{Numerical tests}
\label{sec:results}

The solution method is tested using a function which can be varied to
examine the performance of the algorithm with regard to potential
sources of error. The main sources of error in the algorithm arise
from the interpolation schemes and aliasing. These errors arise in
both the direct method, where the input is specified on the DHT nodes,
and when the input must be interpolated from another mesh onto these
points. 

Interpolation errors arise when the interpolation scheme is unable to
accurately resolve the function which is being interpolated. This can
occur because the interpolation method proper does not have the
required properties to give a well-converged estimate of the
underlying function, or because the interpolation nodes are not dense
enough to take advantage of an otherwise good interpolation method. In
this sense, `interpolation scheme' refers both to the explicitly
stated polynomial interpolation method used to transfer data from the
input mesh to the DHT nodes, and to the interpolation which is
performed implicitly in the quadrature scheme of the DHT. 

The second source of error is aliasing, when the point distribution is
not dense enough to capture the spatial frequencies present in the
input. This can happen in the DHT, if the analytically defined input
has wavenumbers $\alpha_{m}$ with $m>M$, so that the expansion of
Equation~\ref{equ:dht} does not contain the full set of radial
wavenumbers $\alpha_{m}$ present in the input. Clearly, it can also
happen in the Chebyshev interpolation scheme if the node density is
not high enough, even if the DHT would otherwise contain enough nodes
to capture the full behavior of the input.

\begin{figure}
  \centering
  \includegraphics{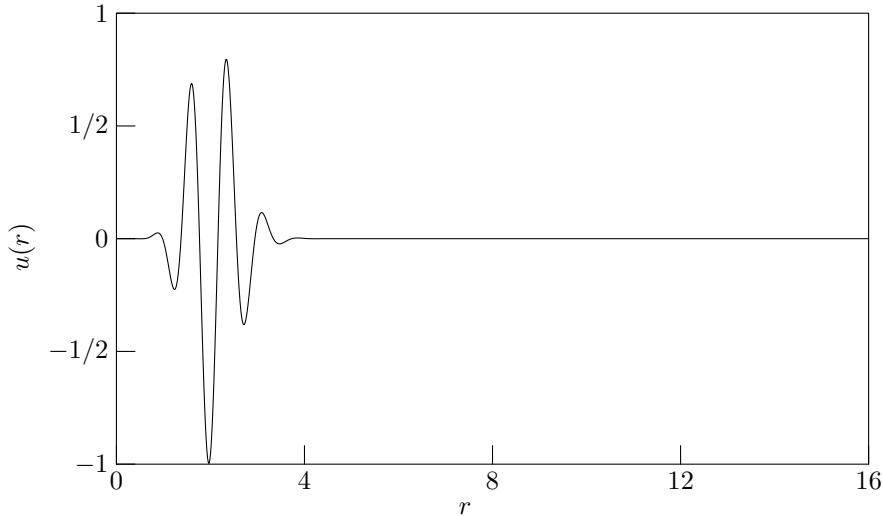}
  \caption{Test function with $\alpha=1$, $\beta=8$, $m=8$}
  \label{fig:test:func}
\end{figure}

In order to assess the performance of the algorithm against these
criteria, the test function:
\begin{align}
  \label{equ:test:func}
  \un(r) &=
  (r/r_{\max})^{m}\E^{-(r^{2}-r_{\max}^{2})/\alpha^{2}}\cos\beta r,\,  m \geq n,\\
  r_{\max} &= \alpha (m/2)^{1/2}, \nonumber
\end{align}
has been adopted. By varying the parameters $\alpha$ and $\beta$, the
function can be varied from monotonically decaying, as in Pataki and
Greengard's test~\cite{pataki-greengard11}, to oscillatory,
Figure~\ref{fig:test:func}. The $r^{m}$ term is required for validity
of the solution and introducing the terms in $r_{\max}$ scales the
amplitude of the cosine on its maximum, so that the maximum amplitude
of $\un$ is one, reducing errors caused by very large values of
$r^{m}$.

In testing the algorithm, $\alpha=1$, $P=16$, and $R=16$. The
parameters varied were $N$ the number of blocks in $r$, $M$ the number
of DHT nodes, $\beta$ the frequency of $\un$ and $\kappa$. The order
$n$ was tested up to~64, with $m=n$ in the evaluation of
$\un(r)$. Tests were conducted for direct solution on the DHT grid, to
evaluate the performance of the underlying method, and with
interpolation from the Chebyshev nodes. The error measure is the
$L_{\infty}$ norm:
\begin{align}
  \label{equ:error}
  \epsilon &= 
  \frac{\max|\un_{c}(r,\kappa)-\un(r)|}{\max|\un(r)|},
\end{align}
where $\un_{c}(r,\kappa)$ is the computed solution. For clarity,
errors greater than~1 have been set to~1 on the plots.

\subsection{Solution on DHT nodes}
\label{sec:solution:dht}

\begin{figure}
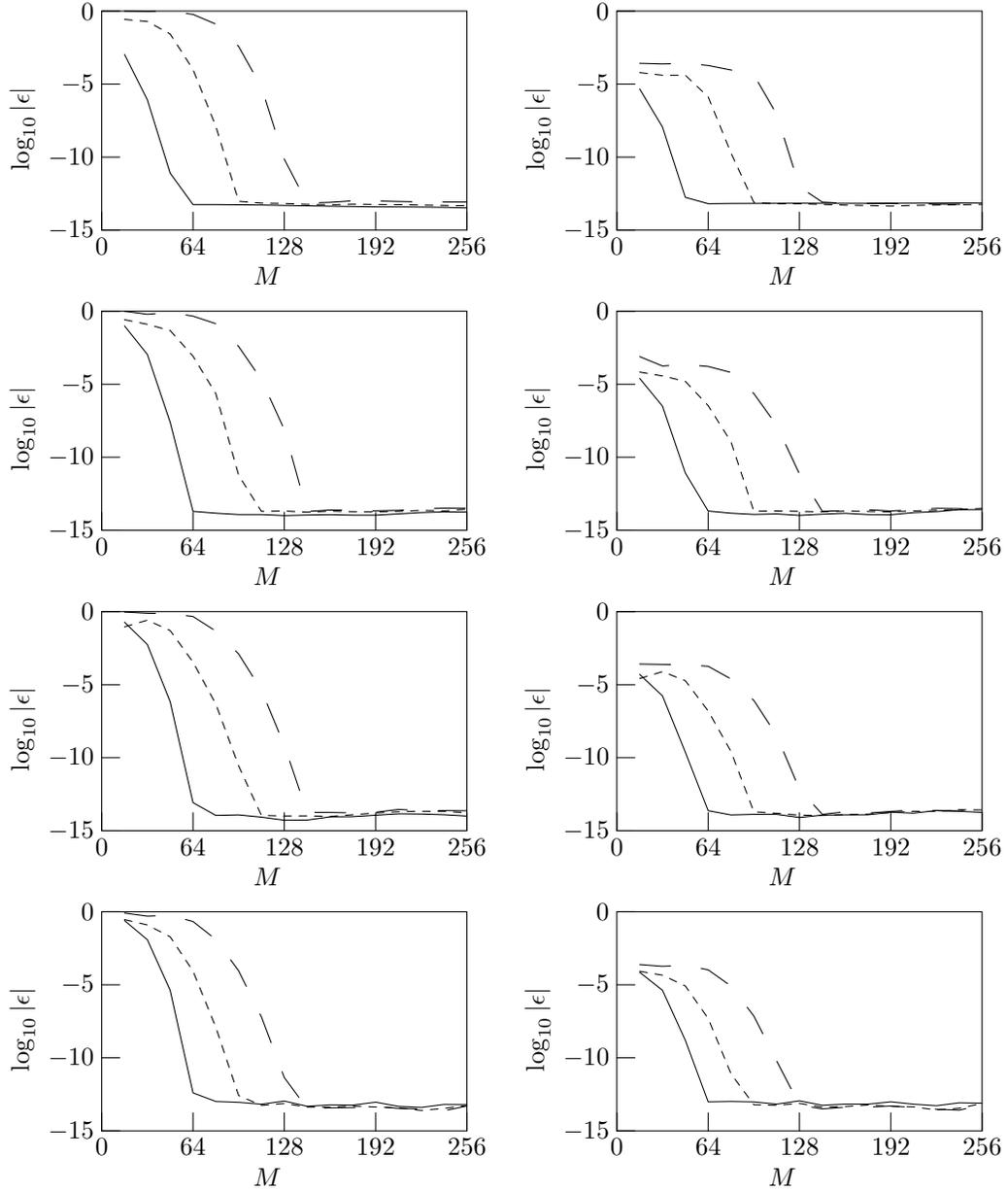

  \centering
  \begin{tabular}{cc}
    \includegraphics{sisc11-figures.1} &
    \includegraphics{sisc11-figures.3} \\
    \includegraphics{sisc11-figures.4} &
    \includegraphics{sisc11-figures.6} \\
    \includegraphics{sisc11-figures.7} &
    \includegraphics{sisc11-figures.9} \\
    \includegraphics{sisc11-figures.10} &
    \includegraphics{sisc11-figures.12}
  \end{tabular}
  \caption{Error in computation on Hankel transform nodes: solid line
    $\beta=0$; dashed line $\beta=8$; long dashed line
    $\beta=16$. Left hand column: $\kappa=16$; right hand column:
    $\kappa=1024$. From top to bottom: $n=0,16,32,64$.}
  \label{fig:error:dht}
\end{figure}

Figure~\ref{fig:error:dht} shows the error in the solution when the
input is specified directly on DHT nodes, given as a function of the
number of nodes $M$, for $\beta=0,8,16$, $\kappa=16,1024$, and
$n=0,16,32,64$. In each plot, the error behavior is quite similar. For
$\beta=0$, the solution is not oscillatory and there is no aliasing
error in the calculation. Thus, the error drops quickly with
increasing $M$, as the node density increases, and reaches a minimum,
machine precision, around $M=64$. In contrast, the error for the
oscillatory solutions, $\beta=8,16$, remains roughly constant for
small node number, before dropping quickly to machine precision. The
initial failure of the error to reduce with $M$ can be ascribed to
aliasing as there are insufficient points to capture the oscillatory
nature of the input, and the reduction in error with $M$ occurs when
there is no longer aliasing and the error is controlled by the
interpolation method.

\subsection{Solution on Chebyshev nodes}
\label{sec:solution:interpolation}

\begin{figure}
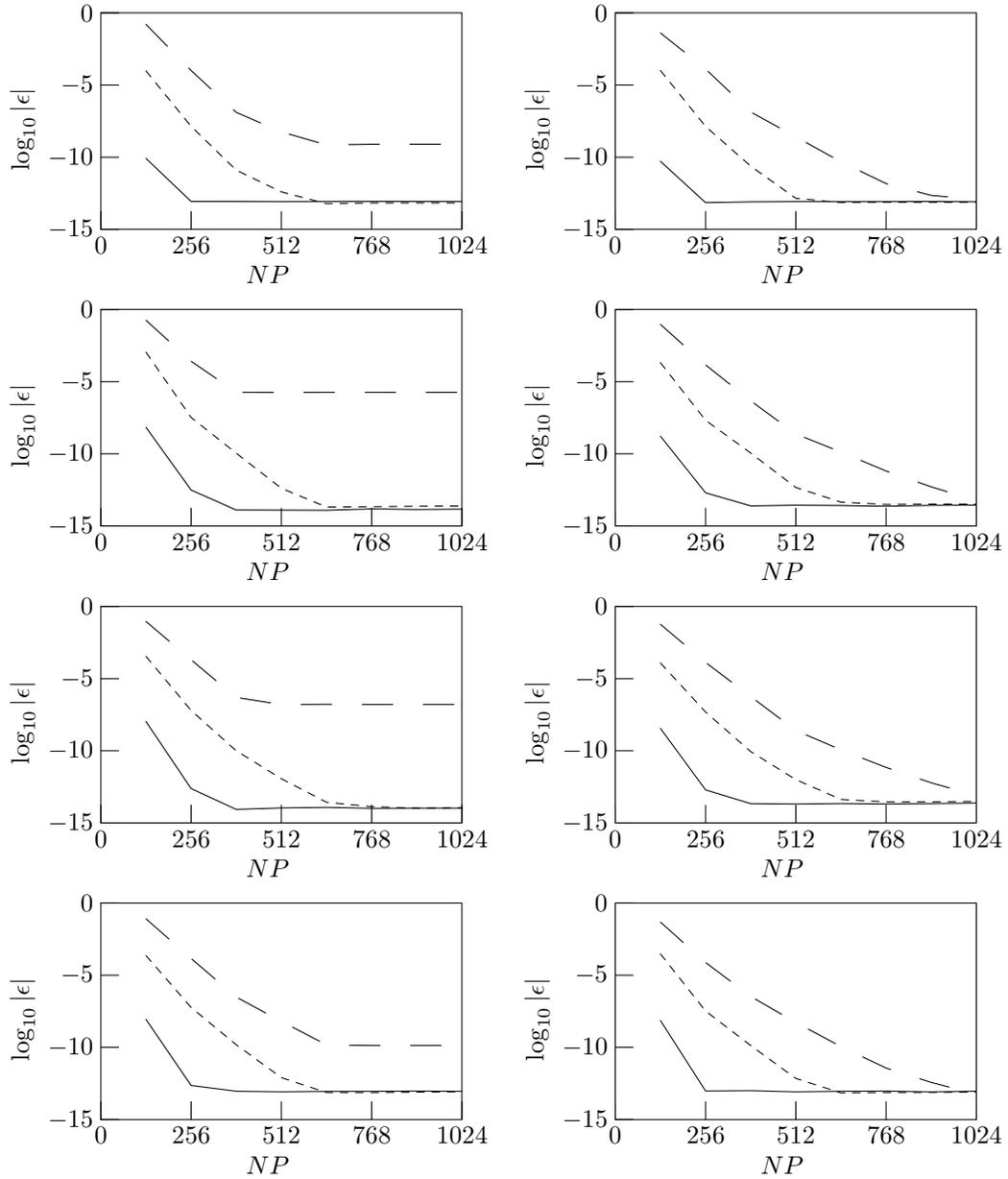

  \centering
  \begin{tabular}{cc}
    \includegraphics{sisc11-figures.103} &
    \includegraphics{sisc11-figures.104} \\
    \includegraphics{sisc11-figures.107} &
    \includegraphics{sisc11-figures.108} \\
    \includegraphics{sisc11-figures.111} &
    \includegraphics{sisc11-figures.112} \\
    \includegraphics{sisc11-figures.115} &
    \includegraphics{sisc11-figures.116}
  \end{tabular}
  \caption{Error in computation on Chebyshev nodes, $\kappa=1024$:
    solid line $\beta=0$; dashed line $\beta=8$; long dashed line
    $\beta=16$. Left hand column: $M=128$; right hand column:
    $M=256$. From top to bottom: $n=0,16,32,64$.}
  \label{fig:error:chebyshev}
\end{figure}

Figure~\ref{fig:error:chebyshev} shows similar data, but for the error
incurred when interpolating from Chebyshev grids onto the DHT
nodes. The parameters varied are the same as in
Figure~\ref{fig:error:dht} but error is now plotted as a function of
the number of nodes $NP$ with $M=128$ and $256$. Referring to
Figure~\ref{fig:error:dht}, the choice $M=128$ gave machine precision
accuracy for $\beta=0$ and $\beta=8$, though not for $\beta=16$, while
$M=256$ gave machine precision errors for all three values of
$\beta$. 

Figure~\ref{fig:error:chebyshev} shows the result of these
choices. The left hand column, where $M=128$ shows the $\beta=0$ and
$\beta=8$ cases reaching the minimum error, though with the error
controlled by $NP$. As the node density increases, the interpolation
scheme reduces the error in the input until the error fixed by the
value of $M$ is reached. For $\beta=16$, the error set by $M$ is
greater than machine precision, and the Chebyshev interpolation scheme
cannot reduce it below the value reached at $NP\approx512$. Referring
to the right hand column, all three cases show steady reduction in
$\epsilon$ down to machine precision, but, as might be expected, the
$\beta=16$ case requires more nodes in order to reduce the error to
the value set by $M$.

\subsection{Computation time}
\label{sec:time}

\begin{figure}
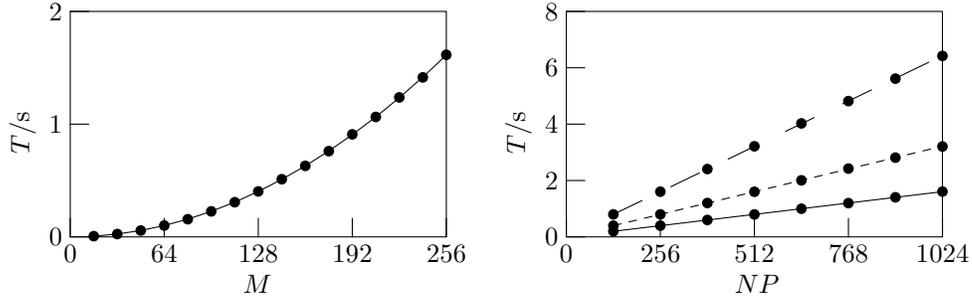

  \centering
  \begin{tabular}{cc}
    \includegraphics{sisc11-figures.1000}
    &
    \includegraphics{sisc11-figures.1100}
  \end{tabular}
  \caption{Computation time for solution with $n=64$, $\kappa=1024$,
    $\beta=16$. Left hand plot uses direct calculation at DHT points;
    circles: computation time; solid line: $M^{2}$ fit. Right hand
    plot shows computation time for Chebyshev nodes; circles:
    computation time; curves: linear fit; solid line: $M=64$; dashed
    line: $M=128$; long dashed line $M=256$.}
  \label{fig:solution:time}
\end{figure}

Finally, Figure~\ref{fig:solution:time} shows the computational time
for the two approaches. The left hand plot shows the time required for
direct solution on the DHT nodes, as a function of $M$. Since the DHT
is implemented as a matrix multiplication, the solution time is
expected to scale as $M^{2}$ and, indeed, the curve fit to the data
points shows a computation time $O(M^{2.00})$. In solving on a
different set of nodes, such as the Chebyshev points used here, the
number of input nodes $NP$ will be greater than $M$, for reasonable
interpolation accuracy. In this case, the computation time can be
expected to scale as $MNP$, i.e. proportional to the number of nodes
and to the number of DHT coefficients used in computing the solution
at each node. This estimate is borne out by the right hand plot in
Figure~\ref{fig:solution:time} where the fits to the computation time
have an exponent equal to~1, to three significant figures, and the
slopes of the lines are proportional to $M$. 

\section{Conclusions}
\label{sec:conclusions}

A method for the solution of a modified Bessel equation which arises
in the solution of Poisson's equation in cylindrical geometries has
been presented, based on the Hankel transform. Numerical testing has
shown the method to be accurate over a wide range of wave numbers and
orders. We conclude that for a proper choice of mesh densities and
Hankel transform order, the method can achieve machine precision
accuracy, for oscillatory and decaying solutions. A sample
implementation of the algorithm, written in C, is available from the
author. 

\section{Acknowledgements}
\label{sec:acknowledgements}

I am grateful to Andras Pataki, Courant Institute, New York
University, for helpful discussions on the implementation of his
Green's function Poisson solver.


\end{document}